%% file: LuzinComb.tex
\documentclass[12pt]{amsart}
\usepackage{amsmath, amssymb, latexsym, amsthm, mathrsfs, color, mathtools, enumitem}
\usepackage[hyphens]{url}
\usepackage[T1]{fontenc}

\usepackage[top=2in, bottom=1.5in, left=1in, right=1in]{geometry}

\input{shortcuts.tex}

\title[Products of Luzin-type sets with combinatorial properties]{Products of Luzin-type sets with combinatorial properties}

\author[P.~Szewczak]{Piotr Szewczak}
\address{Piotr Szewczak, Institute of Mathematics, Faculty of Mathematics and Natural Science College of Sciences, Cardinal Stefan Wyszy\'nski University in Warsaw, W\'oycickiego 1$\slash$3, 01--938 Warsaw, Poland
}
\email{p.szewczak@wp.pl}
\urladdr{http://piotrszewczak.pl}

\author[G.~Wi\'sniewski]{Grzegorz Wi\'sniewski}
\address{Grzegorz Wi\'sniewski, Institute of Mathematics, Faculty of Mathematics and Natural Science College of Sciences, Cardinal Stefan Wyszy\'nski University in Warsaw, W\'oycickiego 1$\slash$3, 01--938 Warsaw, Poland}
\email{grzeorge@o2.pl}

\begin{document}
		
\begin{abstract}
We construct Luzin-type subsets of the real line in all finite powers Rothberger, with a non-Menger product.
To this end, we use a purely combinatorial approach which allows to weaken assumptions used earlier to construct sets with analogous properties.
Our assumptions hold, e.g.,  in the Random model, where already known category theoretic methods fail.
\end{abstract}



\maketitle

\section{Introduction}

\subsection{Combinatorial covering properties}
By \emph{space} we mean a topological space.
A space $X$ is \emph{Rothberger} if, for every sequence  $\eseq{\cU}$ of open covers of the space $X$, there are sets $U_1\in \cU_1, U_2\in \cU_2,\dotsc$ such that the family $\sset{U_n}{n\in\bbN}$  covers $X$.
In the latter case, Rothberger's property can be viewed as a topological version of \emph{strong measure zero}.
We restrict our consideration to \emph{sets of reals}, i.e., spaces homeomorphic to subsets of the real line.
In the realm of sets of reals, we have the following implications.
\[
\text{countability }\longrightarrow\text{ Rothberger }\longrightarrow\text{ strong measure zero}
\]
It is consistent that these properties hold only
for countable spaces and are, thus, equivalent, e.g., in the Laver model for the consistency of Borel's conjecture~\cite{Laver}.
On the other hand, \CH{} implies that the Rothberger property is strictly in between countablility and strong measure zero.

A space $X$ is \emph{Menger} if, for every sequence  $\eseq{\cU}$ of open covers of the space $X$, there are finite subsets $\cF_1\sub \cU_1, \cF_2\sub \cU_2,\dotsc$ such that the family $\bigcup_{n\in\bbN}\cF_n$  covers $X$. 
Menger's property is a generalization of Rothberger's property, and, in the realm of sets of reals, Menger's property is strictly in between $\sigma$-compactness and the Lindel\"of property.

The above covering properties, are central in the \emph{selection principles theory}~\cite{Wiki}.

\subsection{Luzin sets}
A \emph{Luzin} set is an uncountable set of reals whose intersection with every meager set is at most countable; its existence, is independent with ZFC.
A Luzin set is Rothberger (and thus Menger), and it is not $\sigma$-compact.
Thus, a Luzin set is a counterexample to Menger's conjecture which says that, in the realm of sets of reals, the Menger property and $\sigma$-compactness, are equivalent~\cite{Menger24, Hure27}.
A Luzin set also solves Hurewicz's problem, whether the Menger property is equivalent to the Hurewicz property, a formally stronger than Menger, a covering property; every Luzin set is Menger, but not Hurewicz~\cite{Hure27}.
A detailed analysis of Luzin sets, gives a good intuition how to proceed, using combinatorics, in order to construct Menger sets of reals, with no any extra assumptions~\cite{FrMill, BaTs}.

\subsection{Products of Luzin sets}

Let $\cov(\cM)$ be the minimal cardinality of a family of meager subsets of the real line, that covers the real line. 
Let $\cof(\cM)$ be the minimal cardinality of a \emph{cofinal} family of meager subsets of the real line, i.e., every meager subset of the real line is contained in a member of the family. 
A set of reals $X$ is a \emph{$\cov(\cM)$-Luzin} set if $\card{X}\geq\cov(\cM)$ and, for every meager set $M$, we have $\card{X\cap M}<\cov(\cM)$.
Every $\cov(\cM)$-Luzin set is hereditary Rothberger and, assuming that $\cov(\cM)=\cof(\cM)$, a direct construction of such a set is possible.

Luzin sets and $\cov(\cM)$-Luzin sets, play an important role in investigations of products spaces with the Rothberger, or Menger properties.
Just, Miller, Scheepers, and Szeptycki proved that, under \CH{}, there is a Luzin set, whose square is not Menger~\cite[Lemma~2.6]{coc2}; it was the first example showing that products of Rothberger sets, may not be Menger.
Let $\fc$ be the cardinality of the real line.
Assuming \CH{} or the equality $\cov(\cM)=\fc$, independently, many authors proved that there are two $\cov(\cM)$-Luzin sets in all finite powers Rothberger, with a non-Menger product (\cite[Proposition~3.1]{TsAdd}, \cite[Page~205]{coc2}, \cite[Theorem~13]{Sch}, \cite[Section~3]{krawczyk}, \cite[Theorem~4]{bst}).
Category theoretic methods used there, have substantial limitations, which prevent to weaken assumptions.
A natural question arises, whether the existence of such $\cov(\cM)$-Luzin sets can be proven under $\cov(\cM)=\cof(\cM)$.
We show that, assuming that the cardinal number $\cov(\cM)$ is regular, the answer is positive.
In particular, we prove the existence of such sets, e.g., in the Random model, where $\cov(\cM)=\cof(\cM)=\aleph_1<\fc$.
In this model a category theoretic approach fails.

In the second part of the paper, we use already known methods, to construct, for any Luzin-type set $X$, another Luzin type-set $Y$ with some additional properties (as being a topological group, or whose all finite powers are Rothberger) such that the product space $X\x Y$ is not Menger. 

We apply obtained results to products of function spaces with some local properties.
\section{Luzin-type sets via combinatorial approach}

Let $\bbN$ be the set of natural numbers and $\PN$ be the family of subsets of $\bbN$.
Identify every element of the set $\PN$ with its characteristic function, an element of the Cantor cube $\Cantor$.
Using this identification, we introduce a topology in $\PN$.
Since the Cantor cube is homeomorphic to Cantor's set, every subspace of $\PN$ is a set of reals.
Let $\roth$  be the set of all infinite subsets of $\bbN$.
We identify each set $a\in\roth$ with its increasing enumeration, an element of the Baire space $\NN$.
For each natural number $n$, by $a(n)$ denote the $n$-th smallest element of the set $a$.
Thus, we have $\roth\sub\NN$.
Relative topologies in the set $\roth$, from the Baire space $\NN$ and from the Cantor space $\PN$, are equal.
For elements $x,y \in \roth$, we write $x\les y$, if $x$ \emph{is dominated} by $y$, i.e., if the set $\sset{n}{y(n)< x(n)}$ is finite.
A subset of $\roth$  is \emph{dominating} if every element in $\roth$ is dominated by an element from this set.
Analogously, we define relation $\les$ in $\NN$, and a dominating subset of $\NN$.
Let $\fd$ be the minimal cardinality of a dominating set in $\roth$ (equivalently in $\NN$).
For elements $x,y \in \roth$, define $\max\{x,y\} \in \roth$ such that $\max\{x,y\}(n):=\max\{x(n),y(n)\}$, for all natural numbers $n$.
A set $X\sub \roth$ is \emph{2-dominating}, if the set $\sset{\max\{x,y\}}{x,y\in X}$ is dominating.

The cardinal numbers $\cov(\cM)$ and $\cof(\cM)$ do not change if, in their definitions, consider families of meager subsets of $\PN$, instead of the real line.

%Every set of reals of cardinality less than $\cov(\cM)$ is Rothberger, and every set of reals, of cardinality less than $\fd$ is Menger. A space is Menger if no continuous image of this space into $\roth$ (or $\NN$), is dominating~\cite[Proposition~3]{reclaw}.

\begin{thm}\label{thm:main}
Assume that $\cov(\cM)=\cof(\cM)$ and the cardinal number $\cov(\cM)$ is regular.
There are $\cov(\cM)$-Luzin sets $X,Y\sub\roth$ such that:
\be
\item all finite powers of $X$ and $Y$ are Rothberger,
\item the product space $X\x Y$ is not Menger,
\item the union $X\cup Y$ is $2$-dominating.
\ee
\end{thm}

In order to prove Theorem~\ref{thm:main} and later discussion, we need the following auxiliary results and notations.
A set $r\in\roth$ \emph{reaps} a family $A\sub\roth$ if, for each set $a\in A$, both sets $a\cap r$ and $a\sm r$ are infinite.
Let $\fr$ be the minimal cardinality of a family $A\sub\roth$ that no set $r$ reaps.
For elements $x,y\in\roth$, we write $x\leinf y$ if $y\not\les x$, i.e., if the set $\sset{n}{x(n)<y(n)}$ is infinite.
For a set $Z\sub\roth$ and an element $x\in\roth$, we write $Z\leinf x$ if $z\leinf x$ for all elements $z\in Z$.
We use this convention to any binary relation in $\roth$.
For a set $x\in\PN$, define $x\comp:=\bbN\sm x$.
For an element $z\in \roth$ with $z(1)\neq 1$, define $\tz\in\roth$ such that $\tz(1):=z(1)$, and $\tz(n+1):=z(\tz(n))$, for all natural numbers $n$.
For natural numbers $j,k$ with $j<k$, let $[j,k):=\{\,j,j+1,\dotsc,k-1\,\}$. 
For a set $a\in\roth$ and a function $f\in \roth$, define 
\[
a/f:=\sset{n\in\bbN}{a\cap[f(n),f(n+1))\neq\emptyset}.
\]

The following combinatorial characterization of meager subsets of $\PN$ is an important tool, in our constructions.

\blem[{Bartoszy\'nski, Judah~\cite[Theorem~2.2.4]{bartoszynski}}]\label{lem:meager}
Let $M \sub \PN$ be a meager set.
There are a set $a\in\roth$ and a function $f\in \roth$ such that the set $M$ is contained in the meager set
\[
\sset{x\in \PN}{x\cap [f(n),f(n+1))\neq a\cap [f(n),f(n+1)) \text{ for all but finitely many }n}.
\]
\elem

Let $\fd'$ be the minimal cardinality of a set $A\sub \roth$ such that, for each function $b\in \roth$, there is a set $a\in A$ such that, the set  $\sset{n}{\card{a\cap[b(n),b(n+1))}\geq 2}$ is finite.
The following Lemma, used in our main result, is essentially due to Blass~\cite[Theorem~2.10]{blass}.
We provide a proof, for the sake of completeness.

\blem\label{lem:d'}
We have	$\fd=\fd'$.
\elem

\bpf
($\fd\leq \fd'$):
Let $A\sub\roth$ be a set from the definition of the number $\fd'$, of cardinality $\fd'$, closed under finite modifications of its elements.

The set $A$ is dominating in $\roth$:
Fix an element $b\in \roth$.
There is a set $a\in A$ such that the set $\sset{n}{\card{a\cap[b(n),b(n+1))}\geq 2}$ is finite.
Since the set $A$ is closed under finite modifications of its elements, we may assume that, this finite set is empty, and $b(1)<a(1)$.
Fix a natural number $n\in\bbN$, and let $b(n)\leq a(n)$.
If $a(n+1)<b(n+1)$, then $b(n)\leq a(n)<a(n+1)<b(n+1)$, and thus $\card{a\cap[b(n),b(n+1))}\geq 2$, a contradiction.
It follows that $b(n+1)\leq a(n+1)$.
We have $b\les a$, and thus the set $A$ is dominating in $\roth$.

($\fd'\leq \fd$):
Let $D\sub \roth$ be a dominating set, of cardinality $\fd$, closed under finite modifications of its elements.

The set $\sset{\td}{d\in D}$ satisfies the condition from the definition of the number $\fd'$:
Fix an element $b\in \roth$. 
There exists en element $d\in D$ such that $b(n)< d(n)$, for all natural numbers $n$. 
Fix a natural number $n$.
If there is a natural number $k$ with $b(n)\leq \td(k)<\td(k+1)<b(n+1)$, then 
\[
d(\td(k))=\td(k+1)<b(n+1)\leq b(\td(k)),
\]
a contradiction.
Thus, $\sset{n}{\smallcard{\td\cap[b(n),b(n+1))}\geq 2}=\emptyset$.
\epf

We have the following inequalities between considered cardinals~\cite{blass}: 
\[
\cov(\cM)\leq\fd\leq\cof(\cM)\text{ and }\cov(\cM)\leq\fr.
\]

\begin{lem}\label{lem:main}
Let $\cM'$ be a family of meager sets in $\roth$ with $\card{\cM'}<\cov(\cM)$, and $Z\sub\roth$ be a set with $\card{Z}<\cov(\cM)$.
Let $d\in \roth$.
There are elements $x,y\in \roth\sm \bigcup \cM'$ such that $Z\leinf x, Z\leinf y$, and $d\les \max\{x\comp,y\comp\}$.
\end{lem}

\bpf
Let $\cM'=\sset{M_\beta}{\beta < \alpha}$, for some ordinal number $\alpha<\cov(\cM)$.
Fix an ordinal number $\beta<\alpha$.
By Lemma~\ref{lem:meager}, there are a set $a_\beta\in\roth$ and a function $f_\beta\in \roth$ such that the set $M_\beta$ is contained in the meager set
\[
\sset{y\in \roth}{y\cap [f_\beta(n),f_\beta(n+1))\neq a_\beta \cap  [f_\beta(n),f_\beta(n+1)) \text{ for all but finitely many }n}.
\]
Thus, we may assume that the set $M_\beta$ is equal to the above meager set.
Since $\alpha<\cov(\cM)$, there is a set $a\in \roth \sm\bigcup_{\beta<\alpha}M_\beta$.
Then the sets 
\[
f'_\beta:=\sset{f_\beta (n)}{a\cap [f_\beta(n),f_\beta(n+1))=a_\beta \cap [f_\beta(n),f_\beta(n+1)),n\in\bbN}
\]
are infinite, for all ordinal numbers $\beta<\alpha$.

Assume that $Z=\sset{z_\beta}{\beta<\alpha}$.
We have $\alpha <\cov(\cM)=\fd$.
By Lemma~\ref{lem:d'}, there is a function $b\in \roth$ such that the sets
\[
I_\beta:=\sset{n}{\smallcard{\tz_\beta/\td\cap (b(n),b(n+1))}\geq 2},
J_\beta :=\sset{n}{\smallcard{f'_\beta/\td\cap (b(n),b(n+1))}\geq 2}
\]
are infinite for all ordinal numbers $\beta<\alpha$.

Since $\alpha<\cov(\cM)\leq\fr$, there is a decomposition of the set $\bbN$ into infinite sets $r,s,t, u$, each of them reaps the family $\sset{I_\beta,J_\beta}{\beta< \alpha}$. 
Define

\begin{align*}
	x:=\bigcup_{n\in r\cup u}\bigl[\td(b(n)),\td(b(n+1)+1)\bigr)\cup \bigcup_{n\in s} a\cap \bigl[\td(b(n)),\td(b(n+1)+1))\bigr),\\
	y:=\bigcup_{n\in t\cup s}\bigl[\td(b(n)),\td(b(n+1)+1)\bigr)\cup \bigcup_{n\in u} a\cap \bigl[\td(b(n)),\td(b(n+1)+1))\bigr).
\end{align*}

Fix an ordinal number $\beta<\alpha$ and a natural number $n\in s\cap J_\beta$. 
There is a natural number $i$ such that 
\[
f'_\beta(i),f'_\beta(i+1)\in [\td(b(n)),\td(b(n+1)))
\]
Thus,
\[
x\cap [f'_\beta(i),f'_\beta(i+1))=a\cap [f'_\beta(i),f'_\beta(i+1))=a_\beta \cap [f'_\beta(i),f'_\beta(i+1)).
\]
Let $j$ be a natural number such that $f_\beta(j)=f'_\beta (i).$ 
We have
\[
x\cap [f_\beta(j),f_\beta(j+1))=a_\beta \cap [f_\beta(j),f_\beta(j+1)).
\]
Since the set $s\cap J_\beta$ is infinite, the set 
\[
\sset{n\in\bbN}{x\cap [f_\beta(j),f_\beta(j+1))=a_\beta \cap [f_\beta(j),f_\beta(j+1))}
\]
is infinite, too.
Thus, $x\notin\bigcup\cM'.$
Analogously, we get $y\notin\bigcup\cM'$.

Fix an ordinal number $\beta<\alpha$ and a natural number $n\in t\cap I_\beta$. 
Since 
\[
\smallcard{\tz_\beta/\td\cap (b(n),b(n+1))}\geq 2,
\]
there is a natural number $i$ such that \[
\tz_\beta(i),\tz_\beta(i+1)\in [\td(b(n)+1),\td(b(n+1))).
\]
Since $x\cap[\td(b(n)+1),\td(b(n+1)))=\emptyset$, we have
\[
z(\tz_\beta(i))=\tz_\beta(i+1)<\td(b(n+1))\leq x(\tz_\beta(i)),
\]
for all natural numbers $n\in t$.
The set $t\cap I_\beta$ is infinite, and thus $z\leinf x$.
Analogously, we get $z\leinf y$.

Fix a natural number $k\geq \td(b(1))$.
There is a natural number  $n$ such that $k\in [\td(b(n)),\td(b(n+1)))$. 
Assume that $n\in r\cup u$. 
Since $x\comp\cap [\td(b(n)),\td(b(n+1)+1))=\emptyset$, we have 
\[
d(k)<d(\td(b(n+1)))=\td(b(n+1)+1)\leq x\comp(k).
\]
Analogously, if $n\in t\cup s$, then $d(k)<y\comp(k)$.
Thus, $d\les\max\{x\comp,y\comp\}$.
\epf

A \emph{filter} is a subset of $\roth$ with empty intersection, that is closed under finite intersections and taking supersets.
An \emph{ultrafilter} is a maximal filter. 
For elements $x,y\in \roth$, let $[x\leq y]:=\sset{n}{x(n)\leq y(n)}$.
For a set $F\sub \roth$ and elements $x,y\in \roth$, we write $x\leq_F y$, if $[x\leq y]\in F$.
If the above set $F$ is a filter, then the relation $\le_F$ is transitive.
Let $U$ be an ultrafilter.
A set $X\sub \roth$ is \emph{$\leU $-unbounded} if, for each element $b\in\roth$, there is an element $x\in X$ with $b\leU x$.
Let $\fb(U)$ be the minimal cardinality of a $\leU$-unbounded subset of $\roth$.

%\bdfn Let $U$ be an ultrafilter. A set $X\sub\roth$ with $\card{X}\geq \bof(U)$ is a \emph{$U$-scale} if, for each function $b\in\roth$, we have $b\leU x$ for all but less than $\bof(U)$ functions $x\in X$.\edfn

\bdfn[{Szewczak, Tsaban~\cite[Definition~4.1]{ST}}]
Let $U$ be an ultrafilter.
A set $X\sub\roth$ with $\card{X}\geq \bof(U)$ is a \emph{$U$-scale} if, for each element $b\in\roth$, we have 
\[
\card{\sset{x\in X}{x\leU  b}}<\fb(U).
\]
\edfn

\bprp\label{prp:Umain}
Assume that $\cov(\cM)=\cof(\cM)$ and the cardinal number $\cov(\cM)$ is regular.
There are ultrafilters $U, \tU$ and sets $X,Y \sub \roth$ such that:
\be
\item the set $X$ is a $U$-scale and a $\cov(\cM)$-Luzin set,
\item the set $Y$ is a $\tU$-scale and a $\cov(\cM)$-Luzin set,
\item the set $\sset{\max\{x\comp,y\comp\}}{x \in X, y\in Y}$ is dominating in $\roth$.
\ee
\eprp

\bpf
Let $\sset{d_\alpha}{\alpha <\cov(\cM)}$ be a dominating set in $\roth$ and $\sset{M_\alpha}{\alpha <\cov(\cM)}$ be a cofinal family of meager sets in $\PN$.

Let $F_0=\tF_0=\{\bbN\}$ and $x_0,y_0 \in \roth$ be elements such that $d_0\les \max\{{x}_0\comp,{y}_0\comp\}$.
Fix an ordinal number $\alpha<\cov(\cM)$.
We show that there are elements $x_\alpha,y_\alpha \in \roth$ and sets $F_\alpha,\tF_\alpha \sub \roth$ such that:
\be[label=(\roman*)]
\item the sets $F_\alpha,\tF_\alpha$ are closed under finite intersections,
\item $\bigcup_{\beta<\alpha} F_\beta \sub F_\alpha$, $\bigcup_{\beta<\alpha} \tF_\beta \sub \tF_\alpha$,
\item $\card{F_\alpha},\smallcard{\tF_\alpha}<\cov(\cM)$,
\item $\sset{d_\beta,x_\beta}{\beta<\alpha} \leq_{F_\alpha} x_\alpha,\sset{d_\beta,y_\beta}{\beta<\alpha} \leq_{\tF_\alpha} y_\alpha$,
\item\label{eq:Luzin} $x_\alpha,y_\alpha \notin \bigcup_{\beta<\alpha} M_\beta$,
\item\label{eq:dom} $d_\alpha \les \max\{{x}_\alpha\comp,{y}_\alpha\comp\}$:
\ee
Let 
\[
F=\bigcup _{\beta<\alpha} F_\beta\text{ and }\tF=\bigcup _{\beta<\alpha} \tF_\beta.
\]
The cardinal number $\cov(\cM)$ is regular, and thus $\card{F},\smallcard{\tF}<\cov(\cM)$. 
Define
\[
S:=\maxfin\sset{d_\beta,x_\beta}{\beta<\alpha}\text{ and }\tilde{S}:=\maxfin\sset{d_\beta,y_\beta}{\beta<\alpha}.
\]
We have $\smallcard{S},\smallcard{\tilde{S}}<\cov(\cM)$.
For elements $s,f\in\roth$,  let $s\circ f\in\roth$ be an element such that $(s\circ f)(n)=s(f(n))$, for all natural numbers $n$.
By Lemma~\ref{lem:main}, there are elements $x_\alpha, y_\alpha \in \roth$ such that :

\be[label=(\roman*),resume]
\setcounter{enumi}{6}
\item\label{eq:infabove} $\sset{s\circ f}{s\in S\cup \tilde{S}, f\in F\cup \tF} \leinf x_\alpha,y_\alpha$,
\item $x_\alpha,y_\alpha \notin \bigcup_{\beta< \alpha} M_\beta$,
\item $d_\alpha \les \max\{{x}_\alpha\comp,{y}_\alpha\comp\}$.
\ee 

The intersection of finitely many elements of the set $F\cup \smallmedset{[s\leq x_\alpha]}{s\in  S}$ is infinite:
Since the set $F$ is closed under finite intersections, it is enough to show that for an element $f\in F$ and a finite subset $S'\sub S$, the intersection $f\cap\bigcap_{s\in S'}[s\leq x_\alpha]$ is infinite.
Let $f\in F$ and $S'\sub S$ be a finite set.
We have
\[
\bigcap_{s\in S'}
[s\leq x_\alpha]=\bigcap_{s\in S'}\sset{n}{s(n)\leq x_\alpha(n)} \supseteq
\smallmedset{n}{\max_{s\in S'}s(n)\leq x_\alpha(n)}=\bigl[\max[S']\leq x_\alpha\bigr].
\]
By the definition of the set $S$, we have $\max[S']\in S$.
By~\ref{eq:infabove}, we have $\max[S']\circ f\leinf x_\alpha$, and thus
\[
\max[S'](f(n))\leq x_\alpha (n) \leq x_\alpha(f(n)),
\]
for infinitely many natural numbers $n$.
It follows that the intersection $f\cap \bigl[\max[S']\leq x_\alpha\bigr]$ is infinite, and the intersection $f\cap \bigcap_{s\in S'}[s\leq x_\alpha]$ is infinite, too.

Let $F_\alpha$ be the set  $F\cup \sset{[s\leq x_\alpha]}{s\in S}$, closed under finite intersections. 
Since $\bbN \in F$, we have $d_\beta,x_\beta\in S$, for all ordinal numbers $\beta<\alpha$.
Thus, $\sset{d_\beta,x_\beta}{\beta<\alpha} \leq_{F_\alpha} x_\alpha$.
Analogously, define the set $\tF_\alpha$ and show, that $\sset{d_\beta,y_\beta}{\beta<\alpha}\leq_{\tF_\alpha} y_\alpha$.
By the construction, we have $\card{F_\alpha},\smallcard{\tF_\alpha}<\cov(\cM)$.

Let $U$ and $\tilde{U}$ be ultrafilters containing the sets $\bigcup_{\alpha<\cov(\cM)} F_\alpha $ and $\bigcup_{\alpha<\cov(\cM)} \tF_\alpha $, respectively.
Define
\[
X:=\sset{x_\alpha}{\alpha<\cov(\cM)}\text{ and } Y:=\sset{y_\alpha}{\alpha<\cov(\cM)}.
\]

The set $X$ is a $U$-scale: 
By the construction, we have $\card{X}=\cov(\cM)$.
Fix an element $b\in \roth$.
There is an ordinal number $\beta <\cov(\cM)$ such that $b\les d_\beta$.
For all ordinal numbers $\alpha$ with $\beta<\alpha <\cov(\cM)$, we have $d_\beta\leU x_\alpha$.
Thus, $\card{\sset{x\in X}{x\leU  b}}<\cov(\cM).$ 
Analogously, the set $Y$ is a $\tilde{U}$-scale.

The set $X$ is a $\cov(\cM)$-Luzin set:
Let $M\sub\roth$ be a meager set.
There is an ordinal number $\alpha<\cov(\cM)$ such that $M\sub M_\alpha$.
By~\ref{eq:Luzin}, we have
\[
M\cap X\sub M_\alpha\cap X\sub\sset{x_\beta}{\beta<\alpha},
\]
and thus $\card{M\cap X}<\cov(\cM)$.
Analogously, the set $Y$ is a $\cov(\cM)$-Luzin set.

By~\ref{eq:dom}, the set
$\sset{\max \{{x_\alpha}\comp,{y_\alpha}\comp\}}{\alpha<\cov(\cM)}$ is a dominating subset of $\sset{\max \{{x}\comp,{y}\comp\}}{x,y \in X\cup Y}$. 
\epf

Let $\Fin$ be the set of all finite subsets of $\bbN$. 

\begin{proof}[Proof of Theorem~\ref{thm:main}]
Let $X,Y \sub \roth$ be $\cov(\cM)$-Luzin sets from Proposition~\ref{prp:Umain}.
All finite powers of the sets $X\cup\Fin$, and $Y\cup\Fin$ are Rothberger~\cite[Lemma~2.20]{stz}.
Let $\tau\colon\PN\to\PN$ be a homeomorphism such that $\tau(x):=x\comp$.
The sets $\tau[X\cup \Fin]$ and $\tau[Y\cup \Fin]$ are homeomorphic copies of the $\cov(\cM)$-Luzin sets $X\cup \Fin$ and  $Y\cup \Fin$, respectively.
Thus, the sets $\tau[X\cup \Fin]$ and $\tau[Y\cup \Fin]$ are $\cov(\cM)$-Luzin subsets of $\roth$, whose all finite powers are Rothberger.

By Proposition~\ref{prp:Umain}(3), a continuous image of the product space $\tau[X\cup\Fin]\x \tau[Y\cup\Fin]$ into $\roth$, under the map $(z,t) \mapsto \max\{z,t\}$, contains a dominating set $\sset{\max\{x\comp,y\comp\}}{x \in X, y\in Y}$.
Thus, the product space $\tau[X\cup\Fin]\x \tau[Y\cup\Fin]$ is not Menger.

By Proposition~\ref{prp:Umain}(3), the union $\tau[X\cup \Fin]\cup\tau[Y\cup\Fin]$ is a $2$-dominating subset of $\roth$.
\epf

An open family of proper subsets of a space is an \emph{$\w$-cover} of the space if, each finite subset of the space is contained in a member of the family.
A space $X$ is \emph{Scheepers} if, for every sequence of open covers $\eseq{\cU}$ of the space $X$, there are finite subsets $\cF_1\sub \cU_1,\cF_2\sub \cU_2, \dotsc$ such that the family $\sset{\bigcup \cF_n}{n \in \bbN}$ is an $\w$-cover of $X$.
By the result of Tsaban and Zdomskyy, it is consistent with ZFC that, the Menger and Scheepers properties are equivalent~~\cite[Theorem~3.7]{semtrich}.
Assuming that $\fd\leq \fr$, there is a Menger set of reals that is not Scheepers~\cite[Theorem~2.1]{stz}.
A separation of these properties is also possible with Luzin type-sets.
Assuming that the cardinal number $\cov(\cM)$ is regular, the following theorem is a generalization of the result of Bartoszy\'nski, Shelah, and Tsaban~\cite[Theorem~3.1]{bst}.

\bthm\label{thm:bst}
Assume that $\cov(\cM)=\cof(\cM)$ and the cardinal number $\cov(\cM)$ is regular.
There is a $\cov(\cM)$-Luzin set that is not Scheepers.
\ethm

\bpf
Let $X$ and $Y$ be sets from Theorem~\ref{thm:main}. 
The set $X\cup Y$ is a $\cov(\cM)$-Luzin set, a $2$-dominating subset of $\roth$.
No $2$-dominating subset of $\roth$ is Scheepers~\cite[Theorem~2.1]{cbc}.
\epf

\bcor
Assume that $\cov(\cM)=\cof(\cM)$ and the cardinal number $\cov(\cM)$ is regular.
There is a hereditary Rothberger set of reals that is not Scheepers.
\ecor

Assume that $\cov(\cM)=\cof(\cM)$ and the cardinal number $\cov(\cM)$ is regular.
It follows from Proposition~\ref{prp:Umain} that, in particular, there is an ultrafilter $U$ with $\bof(U)=\cov(\cM)$ and a $\cov(\cM)$-Luzin set that is a $U$-scale.
The below Proposition~\ref{prp:ULuzin} generalizes this result.

\bprp\label{prp:ULuzin}
Assume that $\cov(\cM)=\cof(\cM)$.
Let $U\sub \roth$ be an ultrafilter with $\fb(U)=\cov(\cM)$.
There is a $\cov(\cM)$-Luzin set that is a $U$-scale.
\eprp

We need the following Lemma.

\blem\label{lem:ULuzin}
Let $U$ be an ultrafilter with $\fb(U)=\cov(\cM)$.
Let $\cM'$ be a family of meager sets in $\PN$ with $\card{\cM'}<\cov(\cM)$ and $Z\sub \roth$ be a set with $\card{Z}<\cov(\cM)$.
There is an element $x\in \roth\sm \bigcup \cM'$ such that $Z\le_U x$.
\elem

\bpf
Since $\card{Z}<\cov(\cM)=\bof(U)$, there is an element $y\in \roth$ such that $Z\leU y$.
Let $\cM'=\sset{M_\beta}{\beta<\alpha}$ for some ordinal number $\alpha<\cov(\cM)$.
Fix an ordinal number $\beta<\alpha$.
By Lemma~\ref{lem:meager}, there are a set $a_\beta\in\roth$ and a function $f_\beta\in\roth$  such that the set $M_\beta$ is contained in the meager set
\[
\sset{x\in \PN}{x\cap [f_\beta(n),f_\beta(n+1))\neq a_\beta\cap [f_\beta(n),f_\beta(n+1)) \text{ for all but finitely many }n}.
\]
Thus, we may assume that the set $M_\beta$ is equal to the above meager set.
Since $\alpha<\cov(\cM)$, there is a set $a\in\roth\sm\Un\cM'$.
Then the sets
\[
f_\beta':=\sset{f_\beta(n)}{a\cap [f_\beta (n),f_\beta(n+1))=a_\beta\cap [f_\beta (n),f_\beta(n+1))}
\]
are infinite, for all ordinal numbers $\beta<\alpha$.

Since  $\alpha<\cov(\cM)\leq\fd$, by Lemma~\ref{lem:d'}, there is a function $b\in \roth$ such that the sets 
\[
I_\beta:=\sset{n}{\card{f'_\beta\cap [b(n),b(n+1))}\geq 2}
\]
are infinite, for all ordinal numbers $\beta<\alpha$.
We may assume that $y(b(n))<b(n+1)$, for all natural numbers $n$.

Since $\alpha<\cov(\cM)\leq\fr$, there is a decomposition of the set $\bbN$ into infinite sets $r_1,r_2,r_3$, each of them reaps the family $\sset{I_\beta}{\beta\leq\alpha}$.
Since $U$ is an ultrafilter, for the only one of the above sets $r_i$, we have $\bigcup_{n\in r_i} [b(n),b(n+1)) \in U$.
Let say that for $r_1$.
Define
\begin{align*}
r_1^2:=&\sset{n\in r_1}{\min\sset{k}{k\geq n, k\notin r_1} \in r_2},\\
r_1^3:=&\sset{n\in r_1}{\min\sset{k}{k\geq n, k\notin r_1} \in r_3}.   
\end{align*}
We have $r_1=r_1^2\cup r_1^3$ and $r_1^2\cap r_1^3=\emptyset$.
Since $U$ is an ultrafilter, for the only one of the above sets $r_1^j$, we have $\bigcup_{n\in r_1^j} [b(n),b(n+1)) \in U$.
Let say say that for $r_1^2$.
Let
\[
x:=\bigcup_{n\in r_3} [b(n),b(n+1))\cap f.
\]
For a natural number $k\in \bigcup_{n\in r_1^2} [b(n),b(n+1))$, we have
\[
y(k)<y(b(n+1))<b(n+2)\leq x(k).
\]
Then $y\leU  x$, and thus $Z\leU  x$.

Fix an ordinal number $\beta<\alpha$.
Since the set $r_3$ reaps the family $\sset{I_\beta}{\beta\leq \alpha}$, the set $I_\beta\cap r_3$ is infinite.
We have
\[
x\cap[f'_\beta(n),f'_\beta(n+1))=a\cap[f'_\beta(n),f'_\beta(n+1)),
\]
for infinitely many natural numbers $n$.
By the definition of the set $a'_\beta$, we have
\[
x\cap[f_\beta(n),f_\beta(n+1))=a\cap[f_\beta(n),f_\beta(n+1))=a_\beta\cap[f_\beta(n),f_\beta(n+1)),
\]
for infinitely many natural numbers $n$.
Thus, $x\notin M_\beta$.
\epf

\begin{proof}[Proof of Proposition~\ref{prp:ULuzin}]
Let $\sset{M_\alpha}{\alpha < \cov(\cM)}$ be a cofinal family of meager sets in $\PN$ and $\sset{d_\alpha}{\alpha<\cov(\cM)}$ be a dominating set in $\roth$.

Let $x_0 \in \roth$.
Fix an ordinal number $\alpha< \cov(\cM)$. 
By Lemma~\ref{lem:ULuzin}, applied for the family $\sset{M_\beta}{\beta< \alpha}$ of meager sets, and the set $\sset{d_\beta,x_\beta}{\beta<\alpha}$, there is an element $x_\alpha \in \roth\sm\Un_{\beta<\alpha}M_\beta$ such that $\sset{d_\beta,x_\beta}{\beta<\alpha}\leU x_\alpha$.

Let $X:=\sset{x_\alpha}{\alpha <\cov(\cM)}$.
By the construction, we have $\card{X}=\cov(\cM)=\bof(U)$.
Fix an element $b\in \roth$. 
There is an ordinal number $\beta<\cov(\cM)$ such that $b\les d_\beta$.
For every ordinal number $\alpha$ with $\beta<\alpha<\cov(\cM)$, we have $d_\beta\leU  x_\alpha$. 
Then $\card{\sset{x\in X}{x\leU  b}}<\cov(\cM)=\fb(U)$, and thus the set $X$ is a $U$-scale.

By the construction, the set $X$ is a $\cov(\cM)$-Luzin set.
\epf

Let $U$ be an ultrafilter.
A space $X$ is \emph{$U$-Menger}~\cite{ST} if, for every sequence $\eseq{\cU}$ of open covers of the space $X$, there are finite subsets  $\cF_1\sub \cU_1, \cF_2\sub \cU_2,\dotsc$ such that $\sset{n}{x\in \bigcup \cF_n}\in U$, for all points $x\in X$. 

\bcor
Assume that $\cov(\cM)=\cof(\cM)$.
Let $U$ be an ultrafilter with $\fb(U)=\cov(\cM)$.
There is a $U$-Menger $\cov(\cM)$-Luzin set whose all finite powers are Rothberger.
\ecor

\bpf
By Proposition~\ref{prp:ULuzin}, there is a $\cov(\cM)$-Luzin set $X$, that is a $U$-scale.
Then $X\cup\Fin$ is a $U$-Menger $\cov(\cM)$-Luzin set~\cite[Corollary~4.6]{ST}.
By the result of Tsaban, Zdomskyy, and the first named author~\cite[Lemma~2.21]{stz}, all finite powers of this set are Rothberger.
\epf

\section{Luzin-type sets via topological approach}

Assume that $\cov(\cM)=\cf(\fd)$ or $\cov(\cM)=\fd$.
For every $\cov(\cM)$-Luzin set $X$, there is a  Rothberger set of reals $Y$ such that the product space $X\x Y$ is not Menger~\cite[Corollary~2.11]{ST}.
Under various assumptions, we construct, for every $\cov(\cM)$-Luzin set $X$, another $\cov(\cM)$-Luzin set (with some extra properties) such that the product space $X\x Y$ is not Menger.

\subsection{Products of Luzin-type sets in $\ZN$}
Let $\bbZ$ be the set of integers.
Let $f,g\in\ZN$. 
Define functions $(f+g), |f+g|, (f-g) \in\ZN$ as follows:
\begin{align*}
(f+g)(n)&:=f(n)+g(n),\\
|f+g|(n)&:=|f(n)+g(n)|,\\ (f-g)(n)&:=f(n)-g(n),
\end{align*}
for all natural numbers $n$.
For a function $f\in\ZN$ and a set $X\sub\ZN$, let
\[
f-X:=\sset{f-g}{g\in X}.
\]
In the set $\bbZ$, consider the discrete topology.
The space $\ZN$, with the Tychonoff topology, is homeomorphic to the Baire space $\NN$ (and thus also to $\roth$), and $\ZN$ with the group operation $+$, is a topological group.
The cardinal numbers $\cov(\cM)$ and $\cof(\cM)$ do not change if, in their definitions, consider families of meager subsets of $\ZN$, instead of the real line.

\bthm\label{thm:Luzinprod}
Assume that $\cov(\cM)=\cof(\cM)$ and the cardinal number $\cov(\cM)$ is regular.
For every $\cov(\cM)$-Luzin set $X\sub \ZN$, there is  a $\cov(\cM)$-Luzin set $Y\sub \ZN$ such that the product space $X\x Y$ is not Menger.  
\ethm

\bpf
Let $\sset{d_\alpha}{\alpha<\cov(\cM)}$ be a dominating set in  $\NN$ and $\sset{M_\alpha}{\alpha<\cov(\cM)}$ be a cofinal family of meager sets in $\ZN$.
Let
\[
y_0\in d_0-X.
\]
Fix an ordinal number $\alpha<\cov(\cM)$.
The set $d_\alpha - X$ is a $\cov(\cM)$-Luzin set.
Since the set $\bigcup_{\beta< \alpha} M_\beta \cup \sset{y_\beta}{\beta<\alpha}$ is a union of less than $\cov(\cM)$ meager sets, and the cardinal number $\cov(\cM)$ is regular, this union cannot cover the set $d_\alpha - X$.
Thus, there is a function 
\[
y_\alpha \in (d_\alpha- X)\sm \bigl(\bigcup_{\beta< \alpha} M_\beta \cup \sset{y_\beta}{\beta<\alpha}\bigr).
\]

The set $Y:=\sset{y_\alpha}{\alpha < \cov(\cM)}$ is a $\cov(\cM)$-Luzin set: 
By the construction, we have $\card{Y}=\cov(\cM)$.
Let $M\sub \ZN$ be a meager set.
There is an ordinal number $\alpha<\cov(\cM)$ such that $M\sub M_\alpha$.
Since $M\cap Y\sub M_\alpha\cap Y \sub \sset{y_\beta}{\beta <\alpha}$, we have $\card{M\cap Y}<\cov(\cM)$.

The product space $X\x Y$ is not Menger:
For each ordinal number $\alpha<\cov(\cM)$, we have $y_\alpha\in d_\alpha - X$, and thus  there is a function $x_\alpha\in X$ such that $x_\alpha+y_\alpha=d_\alpha$.
Then the continuous image of the product space $X\x Y$ in $\NN$, under the map $(x,y)\mapsto |x+y|$, contains a dominating set $\sset{d_\alpha}{\alpha<\cov(\cM)}$.
Thus, the product space $X\x Y$ is not Menger~\cite[Proposition~3]{reclaw}.
\epf

\bthm\label{thm:Luzinall}
Assume that $\cov(\cM)=\fc$ and the cardinal number $\fc$ is regular.
For every $\fc$-Luzin set $X\sub \ZN$, there is a $\fc$-Luzin set $Y\sub \ZN$, whose all finite powers are Rothberger, and the product space $X\x Y$ is not Menger.
\ethm
\bpf
Let $Q$ be a countable dense subset of $\ZN$.
Let $\ZN=\sset{d_\alpha}{\alpha<\fc}$, and $\sset{M_\alpha}{\alpha<\fc}$ be a cofinal family of meager sets in $\ZN$.
Let $\sset{(\eseq{\cU^\alpha})}{\alpha<\fc}$ be a set of all countable sequences of countable open families in $\ZN$.

By transfinite induction on ordinal numbers $\alpha<\fc$, choose elements $y_\alpha\in\ZN$, sets $Y_\alpha\sub \ZN$ and sets $\eseq{U^\alpha}\sub\ZN$ such that:

\be[label=(\roman*)]
\item $y_\alpha\in (d_\alpha - X)\sm\Un_{\beta<\alpha}M_\beta$,
\item $Y_\alpha:=Q\cup\sset{y_\beta}{\beta<\alpha}$
\item\label{eq:ZNw} If $\eseq{\cU^\alpha}$ is a sequence of $\w$-covers of the set $Y_\alpha$, then the sets $U^\alpha_1\in\cU^\alpha_1, U^\alpha_2\in\cU^\alpha_2,\dotsc$ formulate an $\w$-cover $\sset{U^\alpha_n}{n\in\bbN}$, of the set $Y_\alpha\cup\{y_\alpha\}$.
\item If the condition from~\ref{eq:ZNw} does not hold, then $U^\alpha_n=\ZN$, for all natural numbers $n$.
\item If $\sset{U^\beta_n}{n\in\bbN}\neq\{\ZN\}$, then it is an $\w$-cover of the set $Y_\alpha\cup\{y_\alpha\}$, for all ordinal numbers $\beta<\alpha$.
\ee

Let $Y_0:=Q$.
Fix an ordinal number $\alpha<\fc$, and let $Y_\alpha:=\sset{y_\beta}{\beta<\alpha} \cup Q$.
If the sequence $\eseq{\cU^\alpha}$ is not a sequence of $\w$-covers of the set $Y_\alpha$, then $U^\alpha_n:=\ZN$, for all natural numbers $n$.
Assume that $\eseq{\cU^\alpha}$ is a sequence of $\w$-covers of the set $Y_\alpha$.
Since $\card{Y_\alpha}<\fc=\cov(\cM)$, there are sets $U^\alpha _1 \in \cU^\alpha _1, U^\alpha_2\in\cU^\alpha_2,\dotsc$ such that the family $\sset{U^\alpha_n}{n\in \bbN}$ is an $\w$-cover of the set $Y_\alpha$~\cite[Theorem~4.8]{coc2}.
Assume that, for each ordinal number $\beta<\alpha$, the family $\sset{U^\beta_n}{n\in\bbN}$ is an $\w$-cover of the set $Y_\alpha$, or $\sset{U^\beta_n}{n\in\bbN}=\{\ZN\}$.
For each finite set $F\sub Y_\alpha$ and each ordinal number $\beta<\alpha$, define 
\[
G^{F,\beta}:=\bigcup \set{U^\beta_n}{n \in \bbN, F\sub U^\beta_n}.
\]
The set $G^{F,\beta}$ is open and dense in $\ZN$.
Define
\[
Z_\alpha:=\bigcup_{\beta<\alpha} M_\beta \cup \bigcup\set{\ZN \sm G^{F,\beta}}{\beta\leq \alpha, F\text{ is a finite subset of }X_\alpha}.
\]
The set $Z_\alpha$ is a union of less than $\fc$ meager sets.
Since the cardinal number $\fc$ is regular, the set $Z_\alpha$ cannot cover the Luzin set $d_\alpha- X$, and thus there is a function
\[
y_\alpha \in (d_\alpha-X) \sm Z_\alpha.
\]
Fix an ordinal number $\beta\leq\alpha$.
Assume that $\sset{U^\beta_n}{n\in\bbN}\neq\{\ZN\}$.
The family $\sset{U^\beta_n}{n\in\bbN}$ is an $\w$-cover of the set $Y_\alpha\cup\{y_\alpha\}$:
For every finite set $F\sub Y_\alpha$, we have $y_\alpha \in G^{F,\beta}$.
Thus, there is an ordinal number $n$ with $F\cup \{y_\alpha\} \sub U^\beta _n$.

Let $Y:=\Un_{\alpha<\fc}Y_\alpha$.
By the construction, the set $Y$ is a $\fc$-Luzin set.
For each ordinal number $\alpha<\fc$, we have $y_\alpha \in d_\alpha- X$, and thus there is a function $x_\alpha\in X$ such that $x_\alpha+y_\alpha=d_\alpha$.
Then $\ZN$ is a continuous image of the product space $X\x Y$, under the map $(x,y)\mapsto x+y$.
Since the Menger property is preserved under continuous maps, and the space $\ZN$ is not Menger (it is homeomorphic to $\NN$), the product space $X\x Y$ is not Menger as well.

All finite powers of the set $Y$ are Rothberger:
It is enough to show that, for every sequence $\eseq{\cU}$ of $\w$-covers of the set $Y$, there are sets $U_1\in\cU_1,U_2\in\cU_2,\dotsc$ such that the family $\sset{U_n}{n\in\bbN}$ is an $\w$-cover of $Y$~\cite[Theorem~3.9]{coc2}.
Let $\cU _1,\cU_2,\dotsc $ be a sequence of $\w$-covers of the set $Y$, that are families of open sets in $\ZN$.
There is an ordinal number $\alpha<\fc$ such that the sequence $\eseq{\cU}$ is equal to the sequence $\eseq{\cU^\alpha}$.
By the construction, the sets $U^\alpha_1\in\cU^\alpha_1, U^\alpha_2\in\cU^\alpha_2,\dotsc$ formulates an $\w$-cover $\sset{U^\alpha_n}{n\in \bbN }$ of the set $Y_\beta$, for all ordinal numbers $\beta<\fc$.
Thus, the family $\sset{U^\alpha_n}{n\in \bbN}$ is an $\w$-cover of the set $Y$.
\epf

\subsection{Products of Luzin-type groups in $\PN$} 

A \emph{$\cov(\cM)$-Luzin group}, is a $\cov(\cM)$-Luzin set that is a topological group.

Let $a,b\in \PN$ and $A,B\sub \PN$.
Define
\begin{align*}
a\oplus b:=&  (a\cup b)\sm (a\cap b),\\
 a\oplus B:= &\sset{a\oplus b}{b\in B},\\
  A\oplus B:= &\sset{a\oplus b}{a\in A,b\in B}.
\end{align*}
The space $\PN$, with the group operation $\oplus$, is a topological group.
For $Y\sub\PN$, and $y\in\PN$,
let $\gengp{Y}$ be the group generated by the set $Y$, and $\gengpel{A}{x}:=\gengp{A\cup\{x\}}$.

\bthm\label{thm:Luzingpcomb}
Assume that $\cov(\cM)=\cof(\cM)$.
There exist $\cov(\cM)$-Luzin groups $X,Y\sub \PN$ such that the product space $X\x Y$ is not Menger.
\ethm

\bpf
Let $\sset{d_\alpha}{\alpha<\cov(\cM)}$ be a dominating set in $\roth$.
Let $\sset{M_\alpha}{\alpha < \cov(\cM)}$ be a cofinal family of meager sets in $\PN$, and $\tM_\alpha:=\Un_{\beta<\alpha}M_\beta\cup\sset{x\comp}{x\in\Fin}$, for all ordinal numbers $\alpha<\cov(\cM)$.

By transfinite induction on ordinal numbers $\alpha<\cov(\cM)$, pick elements $x_\alpha$, $y_\alpha\in\PN$, construct increasing groups $X_\alpha,Y_\alpha$ with $\card{X_\alpha},\card{Y_\alpha}<\cov(\cM)$ and groups $\tX_\alpha:=\Un_{\beta<\alpha}X_\beta$, $\tY_\alpha:=\Un_{\beta<\alpha}Y_\beta$ such that

\be[label=(\roman*)]
\item\label{usel2} $x_\alpha,y_\alpha \notin \tM_\alpha\oplus (\tX_\alpha\cup \tY_\alpha)$,
\item\label{usel2a} $d_\alpha \les \max\{{x}_\alpha\comp,{y}_\alpha\comp\}$,
\item\label{emptxy} $(X_\alpha \sm \tX_\alpha) \cap \tM_\alpha = \emptyset,(Y_\alpha \sm \tY_\alpha) \cap \tM_\alpha = \emptyset$:
\ee

Fix an ordinal number $\alpha<\cov(\cM)$.
By Lemma~\ref{lem:main}, applied to the set $\tM_\alpha\oplus (\tX_\alpha\cup \tY_\alpha)$, a union of less than $\cov(\cM)$ meager sets, and to the function $d_\alpha$, there are elements $x_\alpha,y_\alpha\in\roth$ which satisfies~\ref{usel2} and~\ref{usel2a}.
Let $X_\alpha:=\gengpel{\tX_\alpha}{x_\alpha}$ and $Y_\alpha:=\gengpel{\tY_\alpha}{y_\alpha}$.
We have $(X_\alpha \sm \tX_\alpha) \cap \tM_\alpha = \emptyset$:
Suppose not.
Then there are elements $x\in  \tX_\alpha$ and $m \in \tM_\alpha$ such that $x\oplus x_\alpha=m$.
We have $x_\alpha =x \oplus m$, and thus $x_\alpha \in \tX_\alpha \oplus \tM_\alpha$, a contradiction.
Similarly, we have $(Y_\alpha \sm \tY_\alpha) \cap \tM_\alpha = \emptyset$.

By~\ref{emptxy}, the sets $X:=\bigcup_{\alpha < \cov(\cM)} X_\alpha$ and $Y:=\bigcup_{\alpha < \cov(\cM)} Y_\alpha$ are $\cov(\cM)$-Luzin groups.

The product space $X\x Y$ is not Menger:
By~\ref{emptxy} and the definition of the sets $\tM_\alpha$ we have $x\comp\nin \Fin$, for all elements $x\in X\cup Y$.
Thus, the continuous image of the product space $X\x Y$, under the map $(x,y)\mapsto \max\{{x}\comp,{y}\comp\}$, is in $\roth$.
By~\ref{usel2a}, the image contains a dominating set
$\sset{\max \{{x_\alpha}\comp,{y_\alpha}\comp\}}{\alpha<\cov(\cM)}$.
Thus, the product space $X\x Y$ is not Menger~\cite[Proposition~3]{reclaw}.
\epf

\bthm\label{thm:Luzingp}
Assume that $\cov(\cM)=\cof(\cM)$ and the cardinal number $\cov(\cM)$ is regular.
For every $\cov(\cM)$-Luzin group $X\sub \PN$, there is a $\cov(\cM)$-Luzin group $Y\sub \PN$ such that the product space $X\x Y$ is not Menger.
\ethm

\bpf
Let $\sset{d_\alpha}{\alpha<\cov(\cM)}$ be a dominating set in $\roth$ such that $d_\alpha\comp\in\roth$, for all ordinal numbers $\alpha<\cov(\cM)$.
The sets $B_\alpha:=\sset{b \in \roth}{d_\alpha \les b}$ are meager sets of cardinality $\fc$, for all ordinal numbers $\alpha<\cov(\cM)$.
Let $\sset{M_\alpha}{\alpha < \cov(\cM)}$ be a cofinal family of meager sets in $\PN$, and $\tM_\alpha:=\Un_{\beta<\alpha}M_\beta$, for all ordinal numbers $\alpha<\cov(\cM)$.

By transfinite induction on ordinal numbers $\alpha<\cov(\cM)$, pick elements $b_\alpha$, $y_\alpha\in\PN$ and construct increasing groups $Y_\alpha$ with $\card{Y_\alpha}<\cov(\cM)$, and $\tY_\alpha:=\Un_{\beta<\alpha}Y_\beta$ such that:
\be[label=(\roman*)]
\item\label{eq:b} $b_\alpha \in B_\alpha\sm \tY_\alpha\oplus\Fin\oplus X$,
\item $y_\alpha \in (b_\alpha\oplus X)\sm \bigl(\tY_\alpha\oplus \tM_\alpha\cup \sset{y_\alpha}{\beta<\alpha}\bigr)$
\item\label{eq:Luzingp} $(Y_\alpha \sm \tY_\alpha) \cap \tM_\alpha = \emptyset$.
\ee

Since the set $B_0$ is a meager set of cardinality $\fc$, and $\Fin\oplus X$ is a $\cov(\cM)$-Luzin set, there is an element \[
b_0\in B_0\sm (\Fin\oplus X).
\]
Let
\[
y_0\in b_0\oplus X,
\]
and $Y_0:=\gengp{y_0}$.
Fix an ordinal number $\alpha<\cov(\cM)$. 
Let $\tY_\alpha:=\Un_{\beta<\alpha}Y_\beta$.
Since the cardinal number $\cov(\cM)$ is regular, we have $\smallcard{\tY_\alpha}<\cov(\cM)$.
Thus, the set $\tY_\alpha\oplus\Fin\oplus L$, a union of less than $\cov(\cM)$ sets that are $\cov(\cM)$-Luzin sets, is a $\cov(\cM)$-Luzin set.
Since $B_\alpha$ is a meager set of cardinality $\fc$, there is an element
\[
b_\alpha \in B_\alpha\sm (\tY_\alpha\oplus\Fin\oplus X).
\]
Since the set $b_\alpha\oplus X$ is a $\cov(\cM)$-Luzin set and the cardinal number $\cov(\cM)$ is regular, there is an element 
\[
y_\alpha \in (b_\alpha\oplus X)\sm \bigl(\tY_\alpha\oplus \tM_\alpha\cup \sset{x_\beta}{\beta<\alpha}\bigr).
\]
We have $b_\alpha \notin \tY_\alpha\oplus\Fin\oplus X$, the set $X$ is a group, and thus $(y_\alpha\oplus X) \cap (\tY_\beta\oplus\Fin\oplus X)=\emptyset$.
Let $Y_\alpha:= \gengpel{\tY_\alpha}{y_\alpha}$.
 
We have $(Y_\alpha \sm \tY_\alpha) \cap \tM_\alpha = \emptyset$:
Suppose not.
Then there are elements $y\in  \tY_\alpha$ and $m \in \tM_\alpha$ such that $y\oplus y_\alpha=m$.
We have $y_\alpha =y \oplus m$, and thus $y_\alpha \in \tY_\alpha \oplus \tM_\alpha$, a contradiction.
          
By~\ref{eq:Luzingp}, the set $Y:=\bigcup_{\alpha < \cov(\cM)} Y_\alpha$ is a $\cov(\cM)$-Luzin set, and since $\sset{Y_\alpha}{\alpha<\cov(\cM)}$ is an increasing sequence of groups, it is a group.

The product space $X\x Y$ is not Menger:
Let $Z\sub\PN$ be the image of the product space $X\x Y$, under the map $(x,y)\mapsto x\oplus y$.
Suppose that there are elements $x\in X$, $y\in Y$ such that $x\oplus y\in\Fin\sm\{\emptyset\}$.
Then there is an ordinal number $\alpha<\cov(\cM)$ such that $y\in Y_\alpha$.
Thus, there is an element $\ty\in \tY_\alpha$ such that $y=\ty\oplus y_\alpha$.
Since $y_\alpha\in b_\alpha\oplus X$, there is an element $\tx\in X$ such that $y_\alpha=b_\alpha\oplus \tx$.
We have
\[
x\oplus y=x\oplus \ty \oplus y_\alpha =x\oplus \ty \oplus b_\alpha\oplus\tx,
\]
and thus
\[
b_\alpha=\ty \oplus (x\oplus y)\oplus (x\oplus \tx)\in \tY_\alpha\oplus \Fin\oplus X,
\]
a contradiction with~\ref{eq:b}.
Thus, we have $Z\sub\roth\cup\{\emptyset\}$. 
Assume that the product space $X\x Y$ is Menger.
Since the Menger property is preserved under continuous mappings, the set $Z$ is Menger.
Since $\roth \cup \{\emptyset\}$ is a set of reals, the set $Z\sm \{\emptyset\}$ is an $F_\sigma$ subset of $Z$, and thus it is Menger as well.
Since the set $Z\sm\{\emptyset\}$ is a subset of $\roth$ and it contains a dominating set $\sset{b_\alpha}{\alpha <\cov(\cM)}$, it cannot be Menger, a contradiction.
\epf

\section{Applications to function spaces}

For a set of reals $X$, let $\Cp(X)$ be the space of all continuous real-valued functions on the spaces $X$ with the topology of pointwise convergence.
Properties of the set $X$ can describe local properties of the space $\Cp(X)$, and vice versa.
A space $Y$ has \emph{countable fan tightness}~\cite{cft} if, for each point $y\in Y$ and for every sequence $\eseq{U}$ of subsets of the space $Y$ with $y\in\bigcap_n\overline{U_n}$, there are finite sets $F_1\sub U_1, F_2\sub U_2,\dotsc$ such that $y\in\overline{\bigcup_nF_n}$.
If we request that the above sets $\eseq{F}$ are singletons, then the space $Y$ has \emph{countable strong fan tightness}~\cite{sakai}.
A space is \emph{M-separable}~\cite{BBM} if, for every sequence $\eseq{D}$ of dense subsets of the space, there are finite sets $F_1\sub D_1, F_2\sub D_2,\dotsc$ such that the union $\Un_nF_n$ is a dense subset of the space.
If we request that, the above sets $\eseq{D}$ are singletons, then the space is \emph{R-separable}~\cite{BBM}.

Let $X$ be a set of reals.
By the results of Sakai~\cite{sakai},  and Scheepers~\cite[Theorem~13]{coc6}, the statements:
all finite powers of the set $X$ are Rothberger, the space $\Cp(X)$ has countable strong fan tightness, the space $\Cp(X)$ is R-separable, are equivalent.
By the results of Just, Miller, Scheepers, and Szeptycki~\cite[Theorem~3.9]{coc6}, and Scheepers~\cite[Theorem~35]{coc6}, the statements: all finite powers of the set $X$ are Menger, the space $\Cp(X)$ has countable fan tightness, the space $\Cp(X)$ is M-separable, are equivalent, too.

Let $X$ and $Y$ be sets of reals, and $X\sqcup Y$ be a topological sum of these sets.
The product space $X\x Y$ is a closed subspace of the product space $(X\sqcup Y)^2$.
Thus, if $X\x Y$ is not Menger, then $(X\sqcup Y)^2$ is not Menger, too.
Since the product space $\Cp(X)\x \Cp(Y)$ is homeomorphic to the space $\Cp(X\sqcup Y)$, we have the following corollaries from Theorem~\ref{thm:main}.

\bcor
Assume that $\cov(\cM)=\cof(\cM)$ and the cardinal number $\cov(\cM)$ is regular.
There are $\cov(\cM)$-Luzin sets $X$ and $Y$ such that the spaces $\Cp(X)$ and $\Cp(Y)$ have countable strong fan tightness (they are R-separable) but their product space does not have countable fan tightness (is not M-separable).
\ecor 

\bcor
Assume that $\cov(\cM)=\cof(\cM)$ and the cardinal number $\cov(\cM)$ is regular.
There are hereditary Rothberger sets of reals $X$ and $Y$ such that the spaces $\Cp(X)$ and $\Cp(Y)$ have countable strong fan tightness (they are R-separable), but the product space $\Cp(X)\x \Cp(Y)$ does not have countable fan tightness (is not M-separable).
\ecor

\section{Open problems}

In the light of the above results, the following problems arises.

\bprb
Do Theorems~\ref{thm:main},~\ref{thm:bst},~\ref{thm:Luzinprod},~\ref{thm:Luzinall},~\ref{thm:Luzingp} remain true if, we do not assume that the cardinal number $\cov(\cM)$ is regular?
\eprb

\bprb
Assume that $\cov(\cM)=\cof(\cM)$.
Let $X$ be a $\cov(\cM)$-Luzin set.
Is there a $\cov(\cM)$-Luzin group $Y$ such that the product space $X\x Y$ is not Menger?
What if the cardinal number $\cov(\cM)$ is regular?
What if \CH{} holds?
\eprb

\bprb
Let $X$ be a $\cov(\cM)$-Luzin set.
Is there a Menger (Rothberger) set of reals $Y$, whose all finite powers are Menger (Rothberger), but the product space $X\x Y$ is not Menger?
What if \CH{} holds?
\eprb

\bprb
Assume that, the cardinal number $\cov(\cM)$ is singular.
Can a $\cov(\cM)$-Luzin set be a union of less than $\cov(\cM)$ meager sets?
\eprb

\end{document}

%% file: shortcuts.tex
\newcommand{\nc}{\newcommand}

\nc{\thusfar}{\begin{center}{\my{--- Edited thus far ---}}\end{center}}
\nc{\lei}{\le^\oo}
\nc{\card}[1]{\left|#1\right|}
\nc{\medcard}[1]{\biggl|\,#1\,\biggr|}
\nc{\smallcard}[1]{|\,#1\,|}
\nc{\bds}{bidirectional $\roth$-scale}
\nc{\bbN}{\mathbb{N}}%{\w}
\nc{\bbZ}{\mathbb{Z}}
\nc{\beq}{\begin{eqnarray*}}\nc{\eeq}{\end{eqnarray*}}
\nc{\mbq}{\mb{?}}
\nc{\mb}[1]{{\mbox{\textbf{#1}}}}
\nc{\nop}{$\times$}
\nc{\fbn}{\!\!\fbox{\!\nop\!}\!\!}
\nc{\yup}{\checkmark}
\nc{\forces}{\Vdash}
\nc{\name}[1]{\dot{#1}}
\nc{\tf}{\my{FINISHED THUS FAR}}
\nc{\FU}{Fr\'echet--Urysohn}
\nc{\gs}{$\gamma$~space}
\nc{\Ga}{\Gamma}
\nc{\Om}{\Omega}
\nc{\smallbinom}[2]{\begin{psmallmatrix} #1\\ #2 \end{psmallmatrix}}
\nc{\bgamma}{\smallbinom{\Om}{\Ga}}
\newcommand{\two}{\{0,1\}}
\nc{\productive}[2]{\bigl(#1,\allowbreak #2\bigr)^\x}
\nc{\Sel}{\mathsf{S}}
\nc{\sset}[2]{\{\,#1 : #2\,\}}
\nc{\gengpel}[2]{\la#1 , #2\ra}
\nc{\gengp}[1]{\la#1 \ra}
\nc{\smb}[1]{{\!\!\mb{#1}\!\!}}
\nc{\medset}[2]{{\biggl\{\,#1 : #2\,\biggr\}}}
\nc{\smallmedset}[2]{{\bigl\{\,#1 : #2\,\bigr\}}}
\nc{\set}[2]{{\left\{\,#1 : #2\,\right\}}}
\nc{\seq}[2]{{\la\, #1 : #2\,\ra}}
\nc{\eseq}[1]{#1_1, \allowbreak #1_2, \allowbreak\dots} %explicit sequence
\nc{\cube}{(\Cantor)^\bbN}
\nc{\Match}{\op{Match}}
\nc{\concat}[1]{\hat{\phantom{a}}\langle #1\rangle}
\nc{\poset}{\mathbb{P}}
\nc{\fn}[1]{{\op{Fn}(#1\times\w,2)}}
\nc{\linadd}{\op{linadd}}
\nc{\nonprod}{\non^\x}
\nc{\alephes}{{\aleph_0}}
\nc{\my}[1]{\marginpar{\textcolor{red}{***}}\textcolor{red}{#1}}
\nc{\later}[1]{{\color{green} #1}}
\nc{\BTs}[1]{{\color{green} #1 (BT)}}
\nc{\Cp}{\op{C}_\mathrm{p}}
\nc{\Bp}{\op{B}_p}
\nc{\Pa}[8]{\bibitem{#1} {#2}, \emph{#3}, {#4} \textbf{#5} ({#6}), {#7}--{#8}.}
\nc{\tPa}[5]{\bibitem{#1} {#2}, \emph{#3}, {#4}, to appear.}
\nc{\sPa}[4]{\bibitem{#1} {#2}, \emph{#3}, {#4}, submitted.}
\nc{\Bc}[9]{\bibitem{#1} {#2}, \emph{#3}, in: \textbf{#4} (#5), #6 #7, #8--#9.}
\nc{\fD}{\mathfrak{D}}
\nc{\fX}{\mathfrak{X}}
\nc{\Onbd}{\Op_{\mathrm{nbd}}} %{\Op_{\mathsf{nbd}}}
\nc{\Omnb}{\Om_{\mathrm{nbd}}} %{\Om_{\mathsf{nbd}}}
\nc{\od}{\mathfrak{od}}
\nc{\Setting}[7]{\xymatrix@R=4pt@C=7pt{#1\ar@{-}[r]&#2\ar@{-}[r]&#3\\&#4\ar@{-}[u]\\
#5\ar@{-}[uu]\ar@{-}[r] & #6\ar@{-}[u]\ar@{-}[r] & #7\ar@{-}[uu]}}
\nc{\mx}[1]{\begin{matrix}#1\end{matrix}}
\nc{\plim}{p\txt{-}\lim}
\nc{\Bgp}{{\Z^\bbN}}
\nc{\Cgp}{{{\Z_2}^\bbN}}
\nc{\Cite}[1]{\textbf{[#1]}}
\nc{\Next}[1]{{#1^+}}
\nc{\cFin}{\mathrm{cF}}
\nc{\intvl}[2]{{[#1(#2),\allowbreak #1(#2\!+\!1))}}
\nc{\Bdd}{\mathbf{B}}
\nc{\Dfin}{\mathfrak{D}_\mathrm{fin}}
\nc{\grbl}{{\mbox{\textit{\tiny gp}}}}
\nc{\bbP}{\mathbb{P}}
\nc{\BOfat}{\B_{\Om_{\mathrm{fat}}}}%\B_{\mathrm{fat}}}
\nc{\Bgood}{\B_{\mathrm{good}}}
\nc{\compactN}{\cl{\mathbb{N}}}
\nc{\blocks}[2]{\op{cl}_{#2}(#1)}
\nc{\blocksplus}[2]{\op{cl}^+_{#2}(#1)}
\nc{\arx}[1]{\texttt{http://arxiv.org/math/#1}}
\nc{\bq}{\begin{quote}}
\nc{\eq}{\end{quote}}
\nc{\cl}[1]{\overline{#1}}
\nc{\CH}{the Continuum Hypothesis}
\nc{\MA}{Martin's Axiom}
\nc{\Bfat}{\B_\mathrm{fat}}
\nc{\inv}{^{-1}}
\nc{\Cantor}{{\two^\bbN}}%{{2^\w}}
\nc{\bP}{\mathbf{P}}
\nc{\bof}{\op{\fb}}
\nc{\dof}{\op{\fd}}
\nc{\bofF}{\bof(\cF)}
\nc{\sr}[3]{\underset{\mbox{#3}}{\mbox{#1}}}
\nc{\gp}{\binom{\Om}{\Ga}}
\nc{\gpsmall}{\mbox{$\gp$}}
\nc{\gig}{\gimel}%{\gimel\Ga}
\nc{\gns}{\sone(\Om,\gig)}
\nc{\nsr}[2]{#1}
\nc{\Srg}{{\mathbb{S}}}
\nc{\Srgs}{{\mathbb{S}^*}}
\nc{\NN}{{\bbN^{\bbN}}}
\nc{\ZN}{{\Z^{\bbN}}}
\nc{\NNup}{{\bbN^{\uparrow\bbN}}}
\nc{\Pof}{\op{P}}
\nc{\PN}{{\Pof(\bbN)}}
\nc{\roth}{{[\bbN]^{\mbox{\tiny $\infty$}}}} %{{[\w]^{\w}}}
\nc{\Fin}{\mathrm{Fin}}
\nc{\ici}{[\bbN]^{ \infty, \infty}}%{{[\bbN]^{(\aleph_0,\aleph_0)}}}
\nc{\Inc}{{\compactN^{\uparrow\bbN}}}
\nc{\powInc}[1]{{\big(\Inc\big)^{#1}}}
\nc{\powFin}[1]{{\big(\Fin\big)^{#1}}}
\nc{\powPN}[1]{{\big(\PN\big)^{#1}}}
\nc{\NcompactN}{{\compactN^\bbN}}
\nc{\Uarrow}{\smash{\big\uparrow}}
\nc{\LE}{\preccurlyeq}
\nc{\GE}{\succcurlyeq}
\nc{\op}{\operatorname}
\nc{\im}{\op{im}}
\nc{\Span}{\op{span}}
\nc{\maxfin}{\op{maxfin}}
\nc{\ran}{\op{range}}
\nc{\iso}{\cong}
\nc{\Madd}{{\M}^\star}
\nc{\cI}{\mathcal{I}}
\nc{\cJ}{\mathcal{J}}
\nc{\scrA}{\mathscr{A}}
\nc{\scrB}{\mathscr{B}}
\nc{\scrC}{\mathscr{C}}
\nc{\scrD}{\mathscr{D}}
\nc{\scrF}{\mathscr{F}}
\nc{\scrK}{\mathscr{K}}
\nc{\A}{\forall}
\nc{\B}{\mathrm{B}}
\nc{\cB}{\mathcal{B}}
\nc{\bB}{\mathbf{B}}
\nc{\BS}{\mathbf{B}(\mathcal{S})}
\nc{\BF}{\mathbf{B}(\mathcal{F})}
\nc{\BU}{\mathbf{B}(\mathcal{U})}
\nc{\cSp}{\mathcal{S}^+}
\nc{\cFp}{\mathcal{F}^+}
\nc{\cUp}{\mathcal{U}^+}

\nc{\BG}{\B_\Ga}
\nc{\BL}{\B_\Lambda}
\nc{\BT}{\B_\Tau}
\nc{\BTstar}{\B_{\Tau^*}}
\nc{\BO}{\B_\Om}
\nc{\DO}{\cD_\Om}
\nc{\KO}{\cK_\Om}
\nc{\CG}{C_\Ga}
\nc{\CL}{C_\Lambda}
\nc{\CT}{C_\Tau}
\nc{\CTstar}{C_{\Tau^*}}
\nc{\CO}{C_\Om}
\nc{\COgp}{C_{\Om^{\grbl}}}
\nc{\CLgp}{C_{\Lambda^{\grbl}}}
\nc{\BOgp}{\B_{\Om}^{\grbl}}
\nc{\BLgp}{\B_{\Lambda^{\grbl}}}
\nc{\sfC}{\mathsf{C}}
\nc{\sfD}{\mathsf{D}}
\nc{\bD}{\mathbf{D}}
\nc{\Tau}{\mathrm{T}}
\nc{\cA}{\mathcal{A}}
\nc{\cK}{\mathcal{K}}
\nc{\cD}{\mathcal{D}}
\nc{\cF}{\mathcal{F}}
\nc{\cS}{\mathcal{S}}
\nc{\cT}{\mathcal{T}}
\nc{\cG}{\mathcal{G}}
\nc{\cY}{\mathcal{Y}}
\nc{\J}{\mathcal{J}}
\nc{\cL}{\mathcal{L}}
\nc{\cM}{\mathcal{M}}
\nc{\cN}{\mathcal{N}}
\nc{\cH}{\mathcal{H}}
\nc{\cO}{\mathcal{O}}
\nc{\Op}{\mathrm{O}}
\nc{\rmA}{\mathrm{A}}
\nc{\rmF}{\mathrm{F}}
\nc{\rmB}{\mathrm{B}}
\nc{\rmD}{\mathrm{D}}
\nc{\rmP}{\mathrm{P}}
\nc{\cC}{\mathcal{C}}
\nc{\cP}{\mathcal{P}}
\nc{\bbQ}{\mathbb{Q}}
\nc{\bbR}{\mathbb{R}}
\nc{\cU}{\mathcal{U}}
\nc{\Un}{\bigcup}
\nc{\cV}{\mathcal{V}}
\nc{\cW}{\mathcal{W}}
\nc{\Z}{{\mathbb Z}}
\nc{\Impl}{\Rightarrow}
\long\def\forget#1\forgotten{\marginpar{\textcolor{green}{Forgetting...}}}
\nc{\ft}{\mathfrak{t}}
\nc{\fb}{\mathfrak{b}}
\nc{\fc}{\mathfrak{c}}
\nc{\fd}{\mathfrak{d}}
\nc{\fg}{\mathfrak{g}}
\nc{\oo}{\infty}
\nc{\fr}{\mathfrak{r}}
\nc{\fk}{\mathfrak{k}}
\nc{\bidi}{\mathfrak{bidi}}
\nc{\fu}{\mathfrak{u}}
\nc{\fh}{\mathfrak{h}}
\nc{\fp}{\mathfrak{p}}
\nc{\fj}{\mathfrak{j}}
\nc{\fs}{\mathfrak{s}}
\nc{\w}{\omega}
\nc{\W}{\Omega}
\nc{\x}{\times}
\nc{\Iff}{\Leftrightarrow}
\newcommand\comp{^{\text{\tt c}}}
\nc{\nin}{\notin}
\nc{\cat}{\hat{\ }}
\nc{\sub}{\subseteq}
\nc{\spst}{\supseteq}
\nc{\sm}{\setminus}
\nc{\as}{\subseteq^*}%{\let\proclaim\relax}
\nc{\les}{\le^*}
\nc{\leinf}{\le^{\infty}}
\nc{\leS}{\le_S}
\nc{\leF}{\le_{\mathcal{F}}}
\nc{\leU}{\le_{U}}
\nc{\rest}{\restriction}

\nc{\la}{\langle}
\nc{\ra}{\rangle}
\nc{\E}{\exists}
\nc{\dom}{\op{dom}}
\nc{\cov}{\op{cov}}
\nc{\add}{\op{add}}
\nc{\cof}{\op{cof}}
\nc{\cf}{\op{cf}}
\nc{\non}{\op{non}}
\nc{\unif}{\op{non}}
\nc{\COV}{\op{COV}}
\nc{\ADD}{\op{ADD}}
\nc{\COF}{\op{COF}}
\nc{\NON}{\op{NON}}
\nc{\impl}{\to}
\nc{\Lp}{\mathcal{L_\p}}
\nc{\Wlog}{without loss of generality}
\newtheorem{thm}{Theorem}[section]
\nc{\bthm}{\begin{thm}} \nc{\ethm}{\end{thm}}
\newtheorem{prop}[thm]{Proposition}
\nc{\bprp}{\begin{prop}} \nc{\eprp}{\end{prop}}
\newtheorem{fact}[thm]{Fact}
\nc{\bfct}{\begin{fact}} \nc{\efct}{\end{fact}}
\newtheorem{prob}[thm]{Problem}
\nc{\bprb}{\begin{prob}} \nc{\eprb}{\end{prob}}
\newtheorem{lem}[thm]{Lemma}
\nc{\blem}{\begin{lem}} \nc{\elem}{\end{lem}}
\newtheorem{claim}[thm]{Claim}
\nc{\bclm}{\begin{claim}} \nc{\eclm}{\end{claim}}
\newtheorem{cor}[thm]{Corollary}
\nc{\bcor}{\begin{cor}} \nc{\ecor}{\end{cor}}
\newtheorem{conj}[thm]{Conjecture}
\nc{\bcnj}{\begin{conj}} \nc{\ecnj}{\end{conj}}
\theoremstyle{definition}
\newtheorem{defn}[thm]{Definition}
\nc{\bdfn}{\begin{defn}} \nc{\edfn}{\end{defn}}
\newtheorem{obs}[thm]{Observation}
\nc{\bobs}{\begin{obs}} \nc{\eobs}{\end{obs}}
\theoremstyle{remark}
\newtheorem{rem}[thm]{Remark}
\nc{\brem}{\begin{rem}} \nc{\erem}{\end{rem}}
\newtheorem{cnv}[thm]{Convention}
\nc{\bcnv}{\begin{cnv}} \nc{\ecnv}{\end{cnv}}
\newtheorem{exam}[thm]{Example}
\nc{\bexm}{\begin{exam}} \nc{\eexm}{\end{exam}}
\nc{\bpf}{\begin{proof}} \nc{\epf}{\end{proof}}
\nc{\be}{\begin{enumerate}}
\nc{\ee}{\end{enumerate}}
\nc{\bi}{\begin{itemize}}
\nc{\bimy}{\my{\begin{itemize}}
\nc{\eimy}{\end{itemize}}}
\nc{\itm}{\item}
\nc{\ei}{\end{itemize}}
\nc{\Subsection}[1]{\goodbreak\subsection*{#1}}%\ \par}

\nc{\sone}{\mathsf{S}_1}
\nc{\sfin}{\mathsf{S}_\mathrm{fin}}
\nc{\ufin}{\mathsf{U}_\mathrm{fin}}
\nc{\Split}{\mathsf{Split}}

\nc{\gone}{\mathsf{G}_1}    \nc{\gfin}{\mathsf{G}_\mathrm{fin}}

\nc{\tz}{\tilde{z}}
\nc{\tx}{x'}
\nc{\ty}{y'}
\nc{\tU}{\tilde{U}}
\nc{\tH}{\tilde{H}}
\nc{\tM}{\tilde{M}}
\nc{\tY}{\tilde{Y}}
\nc{\tF}{\tilde{F}}
\nc{\td}{\tilde{d}}
\nc{\tX}{\tilde{X}}

\nc{\sep}{
\vspace{2cm}
\noindent
\begin{minipage}{\textwidth}
	\textcolor{red}{\rule{\textwidth}{1pt}}
\end{minipage}
}

%% file: LuzinComb.bbl
\begin{thebibliography}{99}
	
	\bibitem{cft} A.~Arhangel'skii, \emph{Hurewicz spaces, analytic sets and fan-tightness of spaces of functions},
	Soviet Mathematics Doklady \textbf{33 2} (1986), 396--399.
	
	\bibitem{bartoszynski}
	T. Bartoszy\'nski, H. Judah,
	\emph{Set Theory: On the structure of the real line}, 
	A. K. Peters, Massachusetts (1995).
	
	\Pa{bst}{T. Bartoszy\'nski, S. Shelah, B. Tsaban}{Additivity Properties of Topological Diagonalizations}{Journal of Symbolic Logic}{68}{2003}{1254}{1260}
	
	\Pa{BaTs}{T. Bartoszy\'nski, B. Tsaban}{Hereditary topological diagonalizations and the Menger--Hurewicz Conjectures}{Proceedings of the American Mathematical Society}{134}{2006}{605}{615}
	
	\bibitem{BBM}A.~Bella, M.~Bonanzinga, M.~Matveev, \emph{Variations of selective separability}, Topology and its Applications \textbf{156} (2009), 1241--1252.
	
	\bibitem{blass} A.~Blass, \emph{Combinatorial cardinal characteristics of the continuum}, in: \textbf{Handbook of Set Theory} (M. Foreman, A. Kanamori, eds.), Springer, 2010, 395--489.
	
	\Pa{FrMill}{D. Fremlin, A. Miller}{On some properties of Hurewicz, Menger and Rothberger}{Fundamenta Mathematicae}{129}{1988}{17}{33}
	
	\Pa{Canjar}{R. Canjar}{Cofinalities of countable ultraproducts: the existence theorem}{Notre Dame Journal of Formal Logic}{30}{1989}{539}{542}
	
	\Pa{coc2}{W.~Just, A.~Miller, M.~Scheepers, P.~Szeptycki}{The combinatorics of open covers II}{Topology and its Applications}{73}{1996}{241}{266}
	
	\bibitem{Hure27} W. Hurewicz,
	\emph{\"Uber Folgen stetiger Funktionen}, Fundamenta Mathematicae \textbf{9} (1927), 193--204.
	
	\bibitem{krawczyk}A. Krawczyk, H. Michalewski, \emph{Linear metric spaces close to being $\sigma$-compact}, Technical Report 46 (2001) of the Institute of Mathematics, Warsaw University.
	
	\Pa{Laver}{R. Laver}{On the consistency of Borel’s conjecture}{Acta Mathematicae}{137}{1976}{151}{169}
	
	\bibitem{Menger24} K.~Menger,
	\emph{Einige \"Uberdeckungss\"atze der Punktmengenlehre},
	Sitzungsberichte der Wiener Akademie \textbf{133} (1924), 421--444.
	
	\Pa{reclaw}{I. Rec{\l}aw}{Every Luzin set is undetermined in the point-opened game}{Fundamenta Mathematicae}{144}{1994}{43}{54}
	
	\bibitem{sakai} M.~Sakai, \emph{Property C'' and function spaces}, Proceedings of the American Mathematical Society \textbf{104} (1988), 917--919. 
	
	\Pa{Sch}{M. Scheepers}{The length of some diagonalization games}{Archive for Mathematical Logic}{38}{1999}{103}{122}
	
	\bibitem{coc6}M.~Scheepers, \emph{Combinatorics of open covers, VI. Selectors for sequences of dense sets}, Quaestiones Mathematicae \textbf{22} (1999), 109--130.
	
	\Pa{cbc}{M.~Scheepers, B.~Tsaban}{The combinatorics of Borel covers}{Topology and its Applications}{121}{2002}{357}{382}
	
	\Pa{ST}{P. Szewczak, B. Tsaban}{Products of Menger spaces: a combinatorial aproach}{Annals of Pure and Applied Logic}{168}{2017}{1}{18}
	
	\bibitem{stz} P. Szewczak, B. Tsaban, L. Zdomskyy, \emph{Finite powers and products of Menger sets}, submitted, available at \url{http://arxiv.org/abs/1903.03170}.
	
	\bibitem{TsAdd} B.~Tsaban, \emph{Additivity numbers of covering properties}, in: \textbf{Selection Principles and Covering Properties in Topology} (L.~Kocinac, editor), Quaderni di Matematica 18, Seconda Universita di Napoli, Caserta 2006, 245--282.
	
	\bibitem{semtrich} B.~Tsaban, L.~Zdomskyy, \emph{Combinatorial images of sets of reals and semifilter trichotomy}, Journal of Symbolic Logic \textbf{73} (2008), 1278--1288.
	
	\bibitem{Wiki} Wikipedia, \emph{Selection principle},
	\url{https://en.wikipedia.org/wiki/Selection_principle}
	
\end{thebibliography}
